\documentclass[11pt]{amsart}
\usepackage{amsmath}
\usepackage{amsfonts}
\usepackage{graphics}
\usepackage{epstopdf}
\usepackage{epsf}
\usepackage{epsfig}
\usepackage{amssymb}
\usepackage{amscd}
\usepackage{url}

  \title{Tropical Convexity via Cellular Resolutions}
  \author[Florian Block]{Florian Block$^1$}
  \address{Florian Block \newline Technische Universit\"{a}t M\"{u}nchen, Zentrum Mathematik, Boltzmannstr. 3, 85748 Garching, Germany}
  \email{block@in.tum.de}
  \author[Josephine Yu]{Josephine Yu$^2$}
  \address{Josephine Yu \newline Department of Mathematics, University of California,
              Berkeley, CA 94720}
  \email{jyu@math.berkeley.edu}

  \newtheorem{theorem}{Theorem}
  \newtheorem{lemma}[theorem]{Lemma}

   \newtheorem{corollary}[theorem]{Corollary}
  \newtheorem{proposition}[theorem]{Proposition}

\newtheorem{conjecture}[theorem]{Conjecture}\newtheorem{algorithm}[theorem]{Algorithm}

  \theoremstyle{definition}
  
  \newtheorem{example}[theorem]{Example}
   \newtheorem{remark}[theorem]{Remark}

    \newcommand{\N}{{\mathcal N}}

  \renewcommand{\P}{{\mathcal P}}
  \newcommand{\A}{{\mathcal A}}

\newcommand{\TP}{\mathbb{TP}}
  \newcommand{\RR}{{\mathbb R}}
  \newcommand{\ZZ}{{\mathbb Z}}
  \newcommand{\NN}{{\mathbb N}}

\begin{document}

\begin{abstract}
The tropical convex hull of a finite set of points in tropical projective space has a natural structure of a cellular free resolution.  Therefore, methods from computational commutative algebra can be used to compute tropical convex hulls.  Tropical cyclic polytopes are also presented.
\end{abstract}

\date{December 31, 2005}

\thanks{$^1$ This work was carried out while visiting UC Berkeley from fall 2004 to spring 2005.}
\thanks{$^2$ Corresponding author.}


\maketitle

\section{Introduction}

The {\em tropical semiring} $(\RR, \oplus, \odot)$ is the set $\RR$ of real numbers with two binary operations called {\em tropical addition} $\oplus$ and {\em tropical multiplication} $\odot$ defined as
$$
a \oplus b = min(a,b), \text{ and } a \odot b = a+b,  \text{ for all }  
a,b \in \RR.
$$
Then $\RR^n$ has the structure of a semimodule over the semiring $(\RR, \oplus, \odot)$ with tropical addition
$$
(x_1,\dots,x_n) \oplus (y_1,\dots,y_n) = (x_1\oplus y_1, \dots, x_n\oplus y_n),
$$
and tropical scalar multiplication
$$
c \odot (x_1,\dots,x_n) = (c\odot x_1, \dots, c\odot x_n).
$$
A set $A \subset \RR^n$ is called {\em tropically convex} if for all  
$x,y \in A$ and $a,b \in \RR$ also
$
(a\odot x) \oplus (b \odot y) \in A.
$
Notice that we do not put any extra condition on $a$ and $b$ as in usual convexity.
The {\em tropical convex hull} $\text{tconv}(V)$ of a set $V \subset \RR^n$ is the inclusionwise minimal, tropically convex set containing $V$ in $\RR^n$.  Also,
$$
\text{tconv}(V) = \{(a_1 \odot v_1) \oplus \cdots \oplus (a_r \odot v_r) ~:~ v_1, \dots, v_r \in V \text{ and } a_1,\dots,a_r \in \RR\}.
$$

Since any tropically convex set $A$ is closed under tropical scalar multiplication, we identify it with its image under the projection onto the $(n-1)$-dimensional {\em tropical  
projective space}
$$
\TP^{n-1} = \RR^n/(1,\dots,1)\RR.
$$ 
The tropical convex hull of a finite set of points has a natural structure of a polyhedral complex. We refer to \cite{DS} for a more extensive introduction to tropical convexity.

Let $V = \{v_1,\dots,v_r\} \subset \TP^{n-1}$,  $v_i = (v_{i1},\dots,v_{in})$, and $\P=\text{tconv}(V)$.  Let $S=\RR[x_{11},\dots,x_{rn}]$ be the polynomial ring over $\RR$ with indeterminates $x_{ij}$ for $i\in [r] = \{1,\dots,r\}$ and $j\in [n] = \{1,\dots,n\}$.  Let the weight of $x_{ij}$ be $v_{ij}$ and the weight of a monomial ${\bf x^a} = \prod x_{ij}^{a_{ij}} \in S$ be $\sum a_{ij} v_{ij}$. The initial form $\text{in}_V(f)$ of a polynomial $f = \sum c_i {\bf x}^{{\bf a}_i}$ is defined to be the sum of terms $c_i {\bf x}^{{\bf a}_i}$ such that ${\bf x}^{{\bf a}_i}$ has maximal weight.
Let $J$ be the ideal generated by the $2 \times 2$ minors of the  
$r \times n$ matrix $[x_{ij}]$.  Let $I = \text{in}_V (J) = \langle \text{in}_V(f)~:~  
f\in J \rangle$ be the initial ideal of $J$ with respect to $V$.  If $V$ is sufficiently generic, the initial ideal $I$ is a square free monomial ideal.
The square free {\em Alexander dual} $I^*$ of a square free monomial ideal $I = \langle {\bf x}^{{\bf a}^1}, \dots, {\bf x}^{{\bf a}^k} \rangle$ is $$
I^* = m^{{\bf a}^1} \cap \cdots \cap m^{{\bf a}^k},
$$ where each ${\bf a}^i$ is a $0$-$1$ vector and $m^{\bf a} = \langle x_j : a_j = 1 \rangle$.  See \cite{MS, St} for details.
 The following is our main result.

\begin{theorem}\label{thm:main}
For a sufficiently generic set of points $V$ in $\TP^{n-1}$, the  
tropical convex hull $\P = \text{tconv}(V)$ supports a minimal linear  
free resolution of the ideal $I^*$, as a cellular complex.
\end{theorem}

 Moreover, the cellular structure of the minimal free resolution is unique (Remark \ref{rmk:unique}), so we get the following algorithm for computing the  
tropical convex hull of a finite set of points in tropical projective space.

\begin{algorithm}
\label{algo1}
~~
\begin{description}
\item[Input]  A list of points $v_1, \dots, v_r \in \TP^{n-1}$ in  
generic position.
\item[Output]  The tropical convex hull of the input points.
\item[Algorithm]
\begin{itemize}
\item[]
\item[1.] Set $J = \langle 2 \times 2$ minors of the $r \times n$  
matrix $[x_{ij}] \rangle$.
\item[2.] Compute $I = $ {\em in}$_V(J)$.
\item[3.] Compute the Alexander dual $I^*$ of $I$.
\item[4.] Find a minimal free resolution of $I^*$.
\item[5.] Output the desired data about the tropical polytope.
\end{itemize}
\end{description}
\end{algorithm}

\begin{figure}[htbp]
\includegraphics[scale = 0.6]{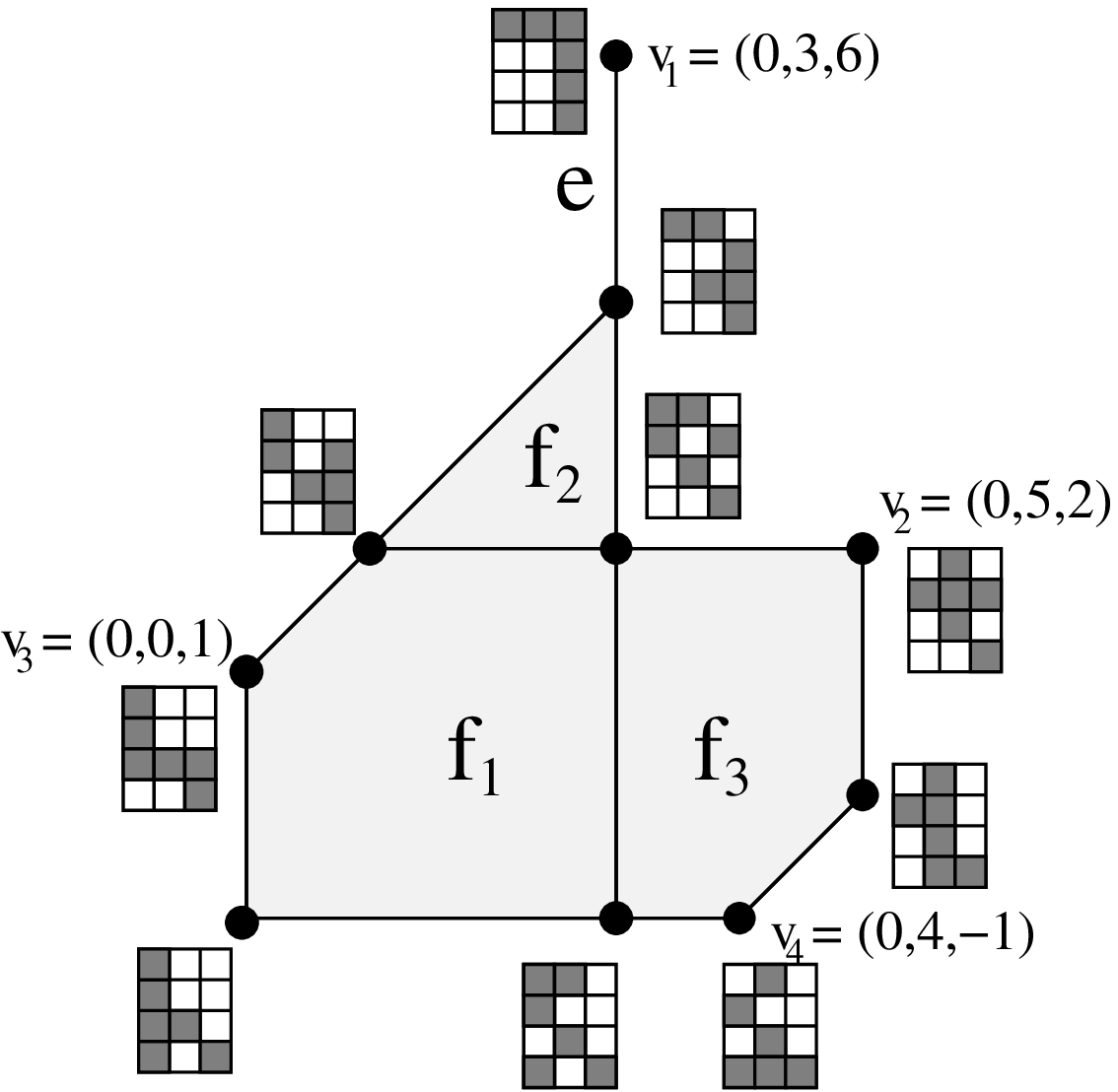}
\caption{Tropical convex hull of four points in $\TP^2$.}
\label{4pt-types}
\end{figure}

A typical output for $r=4$ and $n=3$ is depicted in Figure \ref{4pt-types}.  The ten grids represent ten square free monomials in $S$ of degree six, where each unshaded box represents an indeterminate in the monomial. The cell complex is the minimal free resolution of their ideal.  

Since the set of $2\times 2$ minors of a matrix is fixed under transposition of the matrix, we immediately see the duality between tropical convex hulls of $r$ points in $\TP^{n-1}$ and $n$ points in $\TP^{r-1}$, as shown in \cite{DS}.

The rest of this paper is organized as follows: In Section \ref{sect:proof} we prove Theorem \ref{thm:main} and demonstrate the algorithm with an example. In Section \ref{sect:algo} we deal with algorithmic and computational aspects.  We suggest ways to deal with non-generic points and to get an exterior (halfspace) description of a tropical polytope. We also discuss the efficiency of Algorithm \ref{algo1}. Finally, we study tropical cyclic polytopes in Section \ref{sect:cyclic}.

\section{From Geometry to Algebra and Back}
\label{sect:proof}

We first describe the polyhedral complex structure of the tropical polytope $\P$. Let $W = \RR^{r+n}/(1,\dots,1,-1,\dots,-1)\RR$. Define an unbounded polyhedron as follows:
$$
\P_V = \{(y,z)\in W : y_i + z_j \leq v_{ij} \text{ for all } i \in  
[r], j \in [n]\}.
$$
By \cite{DS}, there is a piecewise linear isomorphism between the  
complex of bounded faces of $\P_V$
and the tropical polytope $\P = \text{tconv}(V)$ given by the projection $(y,z) \mapsto z$.
The boundary complex $\partial \P_V$ of $\P_V$ is  
polar to the regular polyhedral subdivision of the product of simplices  
$\Delta_{r-1}\times\Delta_{n-1}$ induced by the weights $v_{ij}$.  We  
denote this regular subdivision by $(\partial \P_V)^*$.  More  
precisely, a subset of vertices $(e_i,e_j)$ of  
$\Delta_{r-1}\times\Delta_{n-1}$ forms a cell of the subdivision $  
(\partial \P_V)^*$ if and only if the equations $y_i + z_j = v_{ij}$  
indexed by these vertices specify a face of the polyhedron $\P_V$.

Let $\A$ denote the $(r+n) \times rn$ integer matrix whose column vectors are the vertices $(e_i, e_j)$ of $\Delta_{r-1} \times \Delta_{n-1}$, where $i \in [r]$, $j \in  [n]$.  This defines a homomorphism $\ZZ^{rn} \rightarrow  \ZZ^{r+n}$ by $e_{ij} \mapsto  (e_i, e_j)$.
Let $L$ denote its kernel.
The ideal $J$ generated by the $2\times 2$ minors of $[x_{ij}]$ is the (toric) lattice ideal
$$
J = \langle {\bf x^a} - {\bf x^b} ~:~ {\bf a}, {\bf b} \in \NN^{rn}  
\text{ with } {\bf a} - {\bf b} \in L \rangle.
$$
See \cite[Chapter 7]{MS} or \cite[Chapter 8]{St} for details about lattice ideals.

\begin{lemma}
The initial ideal $I$ is independent of the representatives of the points $v_i$ in the tropical projective space.  In other words, if $c\cdot(1,1,\dots,1)$ is added to any $v_i$, the initial ideal $I$ remains the same.
\end{lemma}

\begin{proof}
 The ideal $J$ is homogeneous with respect to any grading assigning the same weight to the variables in each row.
\end{proof}

In the rest of this section, we will assume that the points $v_1,  
\dots, v_r$ are in generic position, i.e., they satisfy the conditions in the next result.

\begin{proposition}
\label{generic_eq}
The following are equivalent.
\begin{itemize}
\item[(1)] The initial ideal $I$ is a monomial ideal.
\item[(2)] The regular subdivision $(\partial P_V)^*$ of $\Delta_{r-1}\times \Delta_{n-1}$ induced by the weights $v_{ij}$ is a triangulation.
\item[(3)] The polyhedron $\P_V$ is simple.
\item[(4)] For any $k$ distinct points in $V$, their projections onto a $k$-dimensional coordinate subspace do not lie in a tropical hyperplane,  for any $2 \leq k \leq n$.
\item[(5)] No $k \times k$ submatrix of the $r \times n$ matrix  
$[v_{ij}]$ is tropically singular, i.e., has vanishing tropical  
determinant (e.g. see \cite{DS}), for any $2 \leq k \leq n$.
\end{itemize}
\end{proposition}

\begin{proof}
\noindent \underline{$(2) \iff (3)$}  follows directly from the polarity between the regular subdivisions of $\Delta_{r-1} \times \Delta_{n-1}$ and $\partial \P_V$. \\
\noindent \underline{$(2) \iff (5)$}  is proven in \cite[Proposition 24]{DS}. \\
\noindent \underline{$(4) \iff (5)$}  is proven in \cite[Lemma 5.1]{RGST}.  
\\
\noindent \underline{$(1) \iff (2)$}:
Statement $(1)$ is equivalent to $V$ being in the interior of a full dimensional cone in the Gr\"{o}bner fan of the lattice ideal $J$.  Statement (2) means that $V$ is in the interior of a full dimenional cone in the secondary fan $\N (\Sigma(\A))$ which is the normal fan of the secondary polytope of $\A$ (for details see \cite{St}).  By \cite[Proposition 8.15(a)]{St}, these two fans coincide if $\A$ is unimodular, i.e., all invertible $rank(\A) \times rank(\A)$ submatrices have the same determinant up to sign.
We will check criterion (iv) of \cite[Theorem 19.3]{Sc} for total unimodularity.  Fix a collection of rows of $\A$. Split it according to containment in the upper $r \times rn$ submatrix of the $(r+n) \times rn$ matrix $\A$. Then the sum of the rows in each part is a $0$-$1$ vector. This implies that all submatrices of $\A$ have determinants $0$ or $\pm 1$, so $\A$ is unimodular.
\end{proof}

It also follows from the unimodularity that all monomial initial ideals of $J$ are square free \cite[Corollary 8.9]{St}.  Let $\Delta_V(J)$ be the {\em initial complex} of $J$, i.e., the simplicial complex whose Stanley-Reisner ideal (see \cite{MS,St}) is $I = \text{in}_V(J)$.  We can identify a square free monomial $m \in S$  with the set of indeterminates $x_{ij}$ dividing $m$.  The vertices of  $\Delta_V(J)$ are $x_{ij}$, and the minimal generators of $I$ are the  minimal non-faces of $\Delta_V(J)$.  Moreover, the minimal generators of the Alexander dual $I^*$ are the complements of the maximal cells of $\Delta_V(J)$.  The following lemma follows immediately from \cite[Theorem 7.33]{MS} or \cite[Theorem 8.3]{St} and establishes a connection between the ideal $J$ and the tropical convex hull.

\begin{lemma}
\label{lem:dual}
We have an isomorphism $\Delta_V(J) \cong (\partial \P_V)^*$, as cell complexes.  In particular, there is  a bijection between maximal cells of $\Delta_V(J)$ and those of $ (\partial \P_V)^*$ induced by $x_{ij} \longleftrightarrow (e_i,e_j)$.
\end{lemma}

We will label the vertices of $\P_V$ by the minimal generators of $I^*$ so that $\P_V$ gives a cellular resolution of $I^*$.  First, we have a general lemma about simple polyhedra which can be proved using \cite[Proposition 4.5]{MS}.

\begin{lemma}[{\cite[Section 4.3.6 and Exercises 4.5-6]{MS}}]
\label{lem:res}
Let $P$ be a simple polyhedron (possibly unbounded) with facets $F_1, \dots, F_m$.  Label each face $G$ of $P$ by ${\bf x}^{{\bf a}_G} = \prod_{F_i \nsupseteq G} x_i \in \RR[x_1, \dots, x_m]$. 
Then the complex of bounded faces of $P$ supports a minimal linear free resolution of the square free monomial ideal generated by the vertex labels. 
\end{lemma}

We will apply this to $\P_V$ to prove Theorem \ref{thm:main} stated in the introduction.

\begin{proof}[Proof of Theorem \ref{thm:main}]
Since $V$ is generic, $\P_V$ is simple.  Hence,  
by Lemma \ref{lem:res}, the tropical convex hull $\P$, which is isomorphic to the complex of bounded faces of $\P_V$, supports a minimal linear free  
resolution of the ideal generated by the monomial labels of its  
vertices.  We only need to show that the labels from Lemma  
\ref{lem:res} coincide with the minimal generators of $I^*$.

The facets $F_{ij}$ of $\P_V$ are defined by equations $y_i + z_j = v_{ij}$.  Let $x_{ij}$ be the indeterminate corresponding  
to $F_{ij}$.  For a square free monomial $m$,
\begin{equation*}
\begin{split}
& m \text{ is a vertex label of } \P_V \\
&\iff  \bigcap\{F_{ij} :  x_{ij} \text{ does not divide } m\}\text{ is a  
vertex of } \P_V\\
&\iff  \{(e_i,e_j) : x_{ij} \text{ does not divide } m\} \text{ is a maximal cell of } (\partial \P_V)^* \\
&\iff  \{x_{ij} : x_{ij} \text{ does not divide } m\}  
\text{ is a maximal cell of } \Delta_V(J) \\
&\iff m \text{ is a minimal generator of } I^*.
\end{split}
\end{equation*}
The third equivalence follows from Lemma \ref{lem:dual}.
\end{proof}

\begin{remark}
\label{rmk:unique}
By construction, the monomial labels are unique, so all the multi-graded Betti numbers are at most one.  This combined with the linearity of the resolution implies that the cellular structure of the minimal free resolution is unique.

However, the multi-graded Betti numbers already determine the tropical polytope because in this case a face $F$ contains a face $G$ if and only if the monomial label of $F$ is divisible by the monomial label of $G$.  Moreover, the vertex labels (the minimal generators of $I^*$) determine all the other monomial labels by Lemma \ref{subface}.
\end{remark}

The {\it  dimension} $dim(U)$ of any subset $U$ of $\TP^{n-1}$ is the affine  
dimension of its projection $\{u \in \RR^n : (0,u) \in U\}$ onto the last $n-1$ coordinates.  

\begin{corollary}
\label{dim} 
For any face $F \subset \P$, $dim(F) = deg({\bf x}^{{\bf a}_F}) - (n-1)(r-1).$
\end{corollary}

\begin{figure}[htbp]
\includegraphics[scale = 0.3]{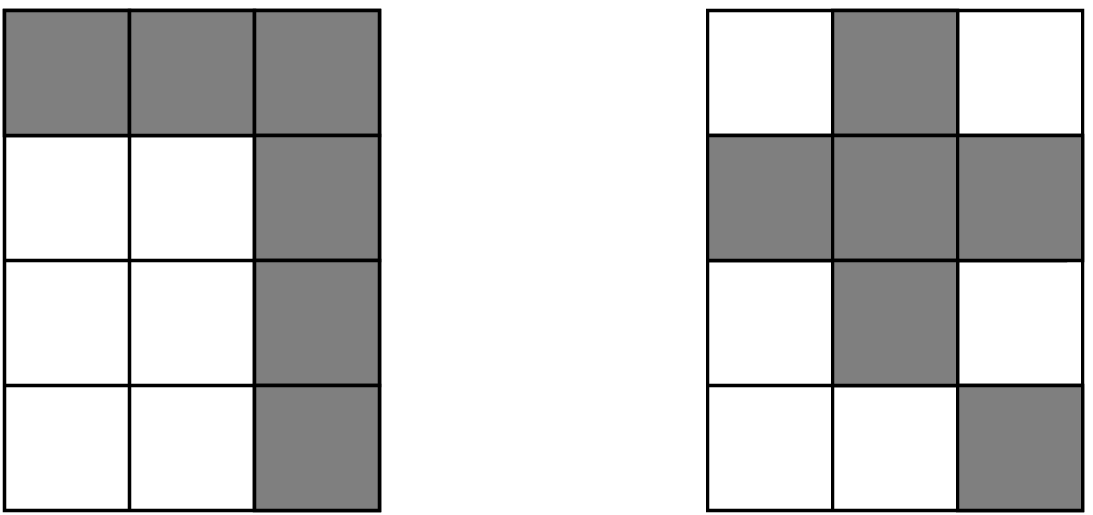}
\caption{Grids representing $x_{21}x_{22}x_{31}x_{32}x_{41}x_{42}$ and $x_{11}x_{13}x_{31}x_{33}x_{41}x_{42}$ for $r=4, n=3$. These are the labels of $v_1$ and $v_2$ in Figure \ref{4pt-types}.}
\label{figure:grid_def}
\end{figure}

The monomial labels have a geometric meaning. 
To have a more intuitive notation, we will represent each squarefree monomial $m \in S$ with an $r \times n$ {\em grid} shaded at position $(i,j)$ if $x_{ij}$ does not divide $m$. Hence the support of ${\bf x}^{{\bf a}_F}$ is left unshaded in the grid (see Figure \ref{figure:grid_def}). 
Let $C_j = cone\{e_i : i\neq j\}$ be the closed cone which is the usual conical (positive) hull of all but one standard unit vector. Suppose $z = (z_1,\dots,z_n) \in \P$ is in the relative interior of a cell with label ${\bf x}^{{\bf a}_z}$, and it is the image of the point $(y,z) \in \P_V$. Then
\begin{equation*}
\begin{split}
x_{ij} \not|~ {\bf x}^{{\bf a}_z} & \iff \left\{
\begin{array}{*{1}{c}}
y_i + z_j = v_{ij}  \\
y_i + z_k \le v_{ik} \quad \forall k
\end{array}
\right\} \iff v_{ij} - z_j \le v_{ik} - z_k \quad \forall k \\
 ~~ & \iff v_i - z \in C_j \iff v_i \in z + C_j.
\end{split}
\end{equation*} 
So the box $(i,j)$ is shaded if and only if the input vertex $v_i$ lies in the {\em sector} $z + C_j$. 
See Figure \ref{fig:hyperplane_halfspace}(b).

This monomial labeling is essentially the same as the labeling by {\em types} introduced in \cite{DS}.  Specifically, for any point $z$ in the relative interior of a cell $F$ in $\P$ with type$(z) = (S_1, \dots, S_n)$, we have $i \in S_j \text{ if and only if } x_{ij}$ does not divide ${\bf x}^{{\bf a}_F}$. The following result follows from \cite[Lemma 10]{DS}.

\begin{lemma}
\label{coord}
Given the monomial label ${\bf x}^{{\bf a}_z}$ of a vertex $z$, its coordinates can be computed by solving the linear system
$$\{z_l - z_k = v_{il} - v_{ik} \; : \; i \in [r],~ k,l \in [n],~ x_{ik} \text{ and } x_{il} \text{ do not divide } {\bf x}^{{\bf a}_z} \}.$$
\end{lemma}

\begin{example}
\label{ex:4pt}
(Four Points in $\TP^2$.)
Assume we are given the following points in $\TP^2$ ($r=4, n=3$):
$$ v_1 =  (0,3,4), v_2 = (0,5,2), v_3 = (0,1,1), v_4 = (0,4,-1).$$
They determine the tropical polytope in Figure \ref{4pt-types}. The points give the weight vector $V =  
(0,3,4,0,5,2,0,1,1,0,4,-1)$ in the polynomial ring $S =  
\RR[x_{11},x_{12},x_{13},x_{21},x_{22},x_{23},x_{31},x_{32},x_{33},x_{41 
},x_{42},x_{43}]$. The initial ideal $I$ and its Alexander dual are
\begin{equation*}
\small
\begin{split}
I = & ~ \text{in}_V\left \langle 2\times 2 \text{ minors of }\left[
\begin{array}{*{3}{c}}
x_{11} & x_{12} & x_{13}  \\
x_{21} & x_{22} & x_{23} \\
x_{31} & x_{32} & x_{33}  \\
x_{41} & x_{42} & x_{43}
\end{array}
\right]\right\rangle \\
= &~  
\langle x_{33}x_{41},x_{23}x_{41},x_{23}x_{31},x_{12}x_{31},x_{31}x_{11}x_{42}, 
x_{11}x_{42},x_{31}x_{41},x_{31}x_{31},x_{13}x_{21}, \\
&  
x_{33}x_{42},x_{22}x_{41},x_{13}x_{32},x_{22}x_{31},x_{11}x_{22},x_{23}x 
_{42},x_{22}x_{33},x_{13}x_{42},x_{13}x_{22},x_{12}x_{21}x_{33}\rangle ,\\
I^* = &~ \langle x_{12}x_{13}x_{22}x_{23}x_{33}x_{42},
x_{12}x_{13}x_{22}x_{23}x_{41}x_{42},
x_{13}x_{22}x_{23}x_{31}x_{33}x_{42},\\
&~~~~ x_{12}x_{13}x_{22}x_{31}x_{41}x_{42},
x_{13}x_{22}x_{31}x_{33}x_{41}x_{42},
x_{13}x_{21}x_{22}x_{31}x_{41}x_{42},\\
&~~~~ x_{11}x_{13}x_{22}x_{23}x_{31}x_{33},
x_{21}x_{22}x_{31}x_{32}x_{41}x_{42},
x_{11}x_{13}x_{31}x_{33}x_{41}x_{42},\\
&~~~~ x_{11}x_{13}x_{23}x_{31}x_{33}x_{41}
 \rangle .
\end{split}
\end{equation*}
Note that $I$ is not generated in degree 2.  Compare the minimal generators of $I^*$ with the grids in Figure \ref{4pt-types}.
The minimal free resolution of $I^*$ is of the form
$$0 \longleftarrow I^* \stackrel{M_0}\longleftarrow S^{10}   
\stackrel{M_1}\longleftarrow S^{12}  \stackrel{M2}\longleftarrow S^3 \longleftarrow  
0.$$ The tropical convex hull  consists of 10 zero-dimensional faces (vertices), 12 one-dimensional  
faces (edges), and 3 two-dimensional faces.

Table \ref{table:M2} shows the monomial matrix $M_2$ in the monomial matrix notation of \cite[Section 1.4]{MS}. The rows correspond to the edges of the tropical polytope, and the three columns, whose labels are omitted here, correspond to the faces $f_1, f_2$, and $f_3$, respectively.

\begin{table}[htbp]
\vskip-0.1in
\begin{equation*}
\scriptsize
\begin{array}{*{2}{c}}
\begin{array}{r}
x_{12}x_{13}x_{22}x_{23}x_{31}x_{33}x_{42} \\
x_{12}x_{13}x_{22}x_{23}x_{33}x_{41}x_{42} \\
x_{12}x_{13}x_{22}x_{23}x_{31}x_{41}x_{42} \\
x_{12}x_{13}x_{22}x_{31}x_{33}x_{41}x_{42} \\
x_{11}x_{13}x_{22}x_{23}x_{31}x_{33}x_{42} \\
x_{12}x_{13}x_{21}x_{22}x_{31}x_{41}x_{42} \\
x_{13}x_{22}x_{23}x_{31}x_{33}x_{41}x_{42} \\
x_{13}x_{21}x_{22}x_{31}x_{32}x_{41}x_{42} \\
x_{11}x_{13}x_{22}x_{31}x_{33}x_{41}x_{42} \\
x_{13}x_{21}x_{22}x_{31}x_{33}x_{41}x_{42} \\
x_{11}x_{13}x_{22}x_{23}x_{31}x_{33}x_{41} \\
x_{11}x_{13}x_{23}x_{31}x_{33}x_{41}x_{42}
\end{array}
&
  \left[
\begin{array}{*{3}{c}}
1 & 0 & 0  \\
-1 & 0 & 0 \\
-1 & 0 & 0  \\
-1 & 1 & 0 \\
0 & 0 & -1 \\
0 & -1 & 0 \\
1 & 0 & -1 \\
{\bf 0} & {\bf 0} & {\bf 0} \\
0 & 0 & 1 \\
0 & 1 & 0 \\
0 & 0 & -1 \\
0 & 0 & 1
\end{array}
\right]
\end{array}
\end{equation*}
\caption{Monomial matrix $M_2$ in Example \ref{ex:4pt}.}
\label{table:M2}
\end{table}
\end{example}

\section{Algorithmic and Computational Aspects}
\label{sect:algo}

Let $0 \leftarrow I^* \leftarrow F_0 \leftarrow \cdots \leftarrow F_m$  
be the free resolution computed by the algorithm, and let $M_i: F_i  
\rightarrow F_{i-1}$ denote the monomial matrices defining the boundary 
maps.
Since the free resolution is linear, the row labels of the 
matrix $M_i$ are in one-to-one correspondence with the faces of dimension  
$i-1$, its column labels with the faces of dimension $i$. An entry in  
$M_i$ is nonzero if and only if its row label divides its column label, which happens if and only if the face corresponding to its column contains the face corrresponding to its row.  Therefore the number of $i$-dimensional faces with $k$ facets in the tropical convex hull is equal to the number of columns of $M_i$ having $k$ nonzero entries.  
 A face $F$ is maximal if and only if it has dimension $n-1$ or the row in $M_{dim(F)+1}$ labeled by ${\bf x}^{{\bf a}_F}$ contains zeroes only. So the eighth row in Table \ref{table:M2} corresponds to edge $e$ in Figure \ref{4pt-types}, which is not contained in any other face.

We can also compute the {\it  
$f$-matrix} $[f_{ij}]$ ($0 \le i \le n-1, 1 \le j$) where $f_{ij}$ is the number of faces having dimension $i$ and $j$ vertices.  We already know the $f$-vector $\sum_j f_{ij}$ which is the sum of columns in the $f$-matrix.
The following result in \cite{DS} was obtained by counting regular triangulations of $\Delta_{r-1} \times \Delta_{n-1}$.
\begin{proposition}[{\cite[Corollary 25]{DS}}]
\label{prop:fvector}
All tropical convex hulls of $r$ generic points in $\TP^{n-1}$ have the same $f$-vector.  The number of faces of dimension $i$ is equal to the multinomial coefficient
$$
{r+n-i-2 \choose r-i-1, n-i-1, i} = \frac{(r+n-i-2)!}{(r-i-1)! \cdot (n-i-1)! \cdot i!}.
$$
\end{proposition}

\subsection*{A Combinatorial Algorithm for Building the Face Poset.}
Given the vertex labels of $\P$, we can compute the whole face poset of $\P$ combinatorially.  The following result follows from \cite[Corollary 14]{DS} and Corollary \ref{dim}.

\begin{lemma}
\label{subface}
Let $F$ be a face of $\P$ with grid $a_F$ and let $b$ be a grid arising from $a_F$ by unshading one box such that no row or column is completely unshaded. Then there is a face $G \supset F$ with label $a_G = b$ and of one dimenstion higher.
\end{lemma}

Conversely, every face can be obtained this way starting from the vertices.  So, instead of computing the free resolution, we can build the fact poset combinatorially if we know the vertex labels, i.e., the minimal generators of $I^*$.  We have an implementation of this algorithm using Macaulay 2 \cite{GS}, Maple \cite{Ma}, and JavaView \cite{Po}.

\begin{remark}[Non-generic Input Vertices]
When the input vertices $V$ are not in generic position, the initial ideal $I$ is not monomial.  In that case, we can replace the weights $V$ with any refinement which makes $I$ a monomial ideal and proceed as before to build the face poset.  We can then compute the coordinates of the vertices using Lemma \ref{coord} and identify vertices with the same coordinates. We suggest this algorithm without having a proof.
\end{remark}

\subsection*{Tropical Halfspaces}

Tropical halfspaces introduced in \cite{Jo} give us an exterior  
description of tropical polytopes.  We can extend our algorithm to find  
such a description.

\begin{figure}[htbp]
\begin{tabular}{ccccccc}
\includegraphics[scale = 0.6]{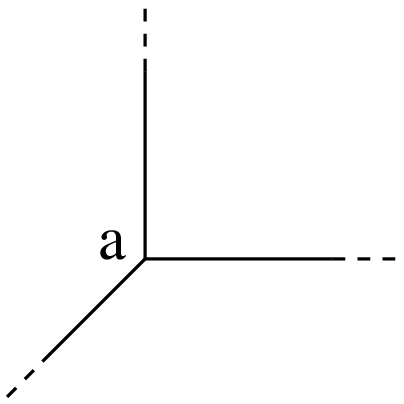} & ~ & ~~&
\includegraphics[scale = 0.6]{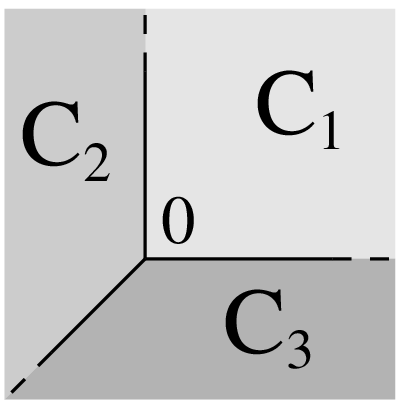} & ~ & ~~ &
\includegraphics[scale = 0.6]{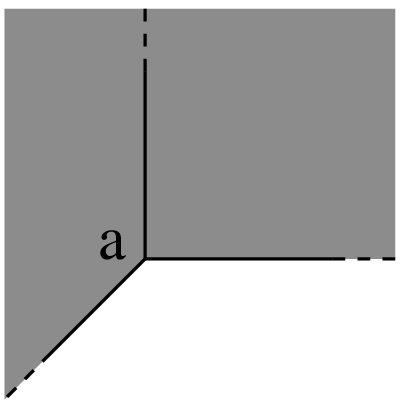}\\
(a) & ~ & ~ & (b) & ~ & ~ & (c)\\
\end{tabular}
\caption{(a) Tropical hyperplane in $\TP^2$ with apex $a$. \newline (b) The sectors at apex $0$ in $\TP^2$. (c) Tropical halfspace $(a,\{1,2\})$ in $\TP^2$. }
\label{fig:hyperplane_halfspace}
\end{figure}

The tropical hyperplane at the apex $a\in \TP^{n-1}$ is the  
set which is the union of boundaries of the {\em sectors} $a + C_i$ (see Figure \ref{fig:hyperplane_halfspace}).
For $a \in \TP^{n-1}$, $\emptyset \neq A \subsetneq [n]$, the set $a + \bigcup_{i \in A} C_i$ is a {\em closed tropical halfspace} (see Figure \ref{fig:hyperplane_halfspace}(c)).
Tropical halfspaces are tropically convex, and a tropical polytope $\P$ is the intersection of the inclusionwise minimal halfspaces containing it \cite{Jo}. 
The apex of such a minimal halfspace must be a vertex of $\P$ on the boundary \cite[Lemma 3.6]{Jo}.
Recall that the box $(i,j)$ in the grid label of a vertex $v$ is shaded if and only if $v_i \in v+C_j$.  Hence  $\P$ is the intersection of the halfspaces $v + \bigcup_{i \in A} C_i$ such that $v$ is a vertex of $\P$ and $A$ is a minimal  subset of columns in the corresponding grid of $v$ such that the shaded boxes in those column cover all the rows. 
This description is redundant in general.  We may be able to refine this result as follows.

\begin{conjecture}
In the generic case, a minimal half space with respect to $\P$ has the form $v + \bigcup_{i \in A} C_i$ where $v$ is a vertex of $\P$ and in the grid label of $v$ the shaded boxes in the columns in $A$ form a partition of $[r]$.
\end{conjecture}

The converse of the conjecture above is not true, i.e., there are  
non-minimal halfspaces of the form described.

\subsection*{Experiments with Computation Time}

We experimented with computing tropical cyclic polytopes $C_{r,n}$ (which will be defined in the next section) with $r$ input vertices in $n-1$ (projective) dimensions.  We used Macaulay 2 \cite{GS} on a Sun Blade 150 (UltraSPARC-IIe 550MHz) computer with 512MB memory.  The computation became infeasible when $r n > 80$ or so, although $r = 30, n = 3$ worked.  The main problem was the insufficient amount of memory.    Some sample computation times for tropical cyclic polytopes are given in Table \ref{table:times}.  We see from the data that computing the Alexander dual can be a problem.  This can be made faster using the {\em Monos Language for Monomial Decompositions} \cite{Mi}.

\begin{table}[htbp]
\begin{tabular}{|c|c|c|c|c|}
\hline
~n~ & ~r~ & Initial ideal & Alexander dual & Free resolution \\
\hline
\hline
3 & 30 & 74 & 433 & 2 \\
4 & 21 & 64 & 944 & 23 \\
6 & 10 & 15 & 221 & 27 \\
8 & 10 & 70 & 4169 & 1106 \\
\hline
\end{tabular}
\caption{Computation times (in seconds) for $I$, $I^*$, and the free resolution for tropical cyclic polytopes $C_{r,n}$.}
\label{table:times}
\vskip-0.2in
\end{table} 

\section{Tropical Cyclic Polytopes}
\label{sect:cyclic}

\begin{figure}[htbp]
\begin{tabular}{ccccc}
\includegraphics[scale = 0.6, viewport=50 70 310 310,clip]{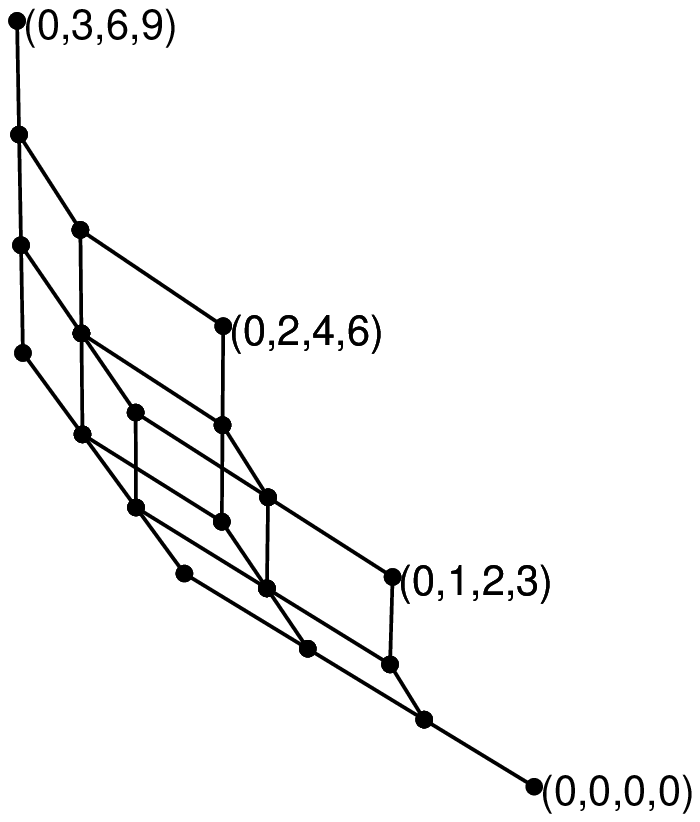}& ~~ & ~ & ~ &
\includegraphics[scale = 0.5]{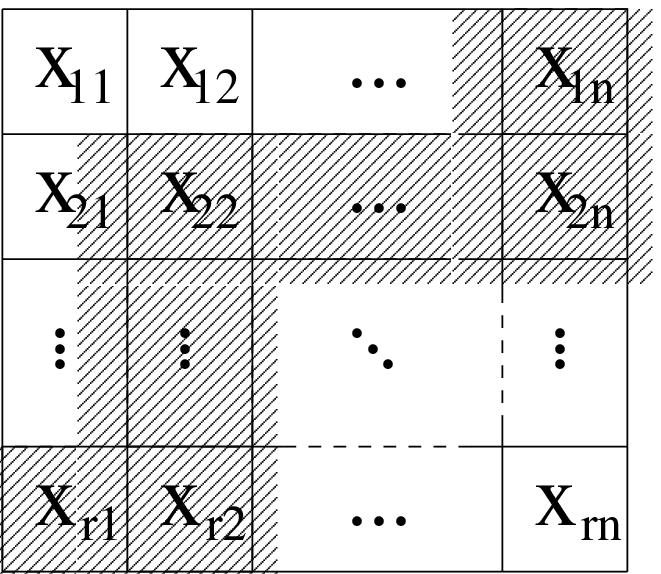}\\
(a) & ~ & ~&~& (b) \\
\end{tabular}
\caption{(a) Tropical cyclic polytope on four vertices in $\TP^3$.
(b) A path corresponding to a generator of the Alexander dual $I^*$.}
\label{fig:cyclic_matrix}
\end{figure}

Define {\em tropical cyclic polytopes} as $C_{r,n} = \text{tconv}\{v_1,  
\dots, v_r\} \subset \TP^{n-1}$,
where $v_{ij} = (i-1)(j-1)$ for $i \in [r],~ j \in [n]$.  Since $(i-1)^{\odot (j-1)} = (i-1)(j-1)$, this is tropical exponentiation.  The $C_{r,n}$ are generic because the minimum in any $k \times k$ minor of the  
matrix $[v_{ij}]$ is attained uniquely by the antidiagonal. An example of a tropical cyclic polytope is shown in Figure \ref{fig:cyclic_matrix}(a).

The $2 \times 2$ minors  of $[x_{ij}]$ form a Gr\"obner basis with respect to $V$, and the initial ideal $I$ is the {\em diagonal initial ideal} generated by the binomials which are on the diagonals of the $2 \times 2$ minors.  This correspond to the {\em staircase triangulation} of $\Delta_{r-1} \times \Delta_{n-1}$.

Consider a path in an $r \times n$ grid representing indeterminates $x_{ij}$, which goes from the lower left corner to the  upper right corner, only moving either right or up at each step as in Figure \ref{fig:cyclic_matrix}(b).  Such paths are precisely the maximal sets, with respect to inclusion, that do not contain diagonal pairs.  Hence their  
complements correspond to the minimal generators of the Alexander dual $I^*$, which are the monomial labels of the vertices of $C_{r,n}$.

\begin{figure}[htbp]
\begin{tabular}{cccc}
\begin{tabular}{ccc}
\includegraphics[scale = 0.5]{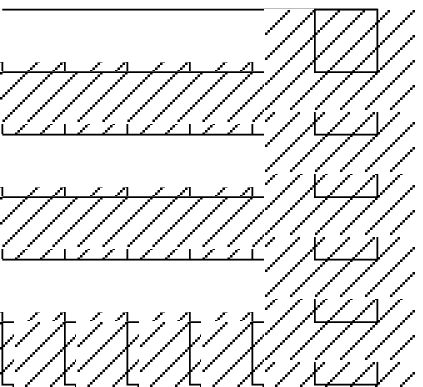} & ~~ &
\includegraphics[scale = 0.5]{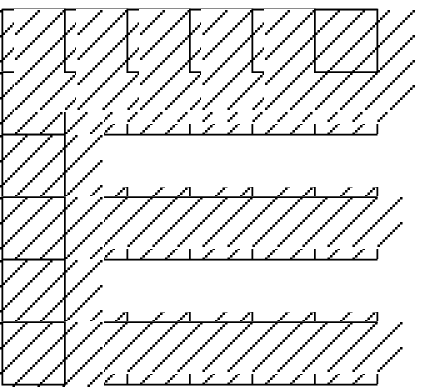}
\end{tabular}& ~~ & ~~ &
\begin{tabular}{c}
\includegraphics[scale = 0.5]{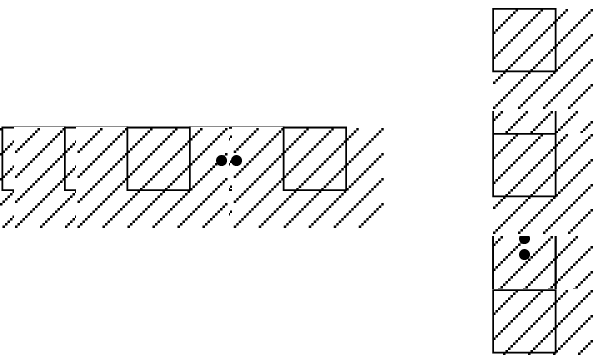}
\end{tabular}\\
(a) & ~~&~~&(b) \\
\end{tabular}
\caption{(a) Paths in grids corresponding to two 1-valent vertices in $C_{r,n}$.  (b) Horizontal and vertical stripes.}
\label{fig:one_valent}
\end{figure}

The labels of the faces of the tropical cyclic polytope  
$C_{r,n}$ are obtained by unshading the boxes on the paths so that the  
remaining shaded set still intersects every row and every column.  For example,  
there are two $1$-valent vertices with grids corresponding to the paths in Figure \ref{fig:one_valent}(a).  The two edges containing these vertices are the only maximal $1$-faces,  
whose labels are obtained by unshading the lower right corner and the  
upper left corner, respectively.

We can identify a vertex of $C_{r,n}$ with the Young diagram above (or  
below) the corresponding path in the $r \times n$ grid.  Then the  
$1$-skeleton of $C_{r,n}$ is the Hasse diagram of the Young lattice of  
the Young diagram fitting in an $(r-1) \times (n-1)$ grid.

\begin{figure}[tbp]
\begin{tabular}{ccc}
\includegraphics[scale = 0.4]{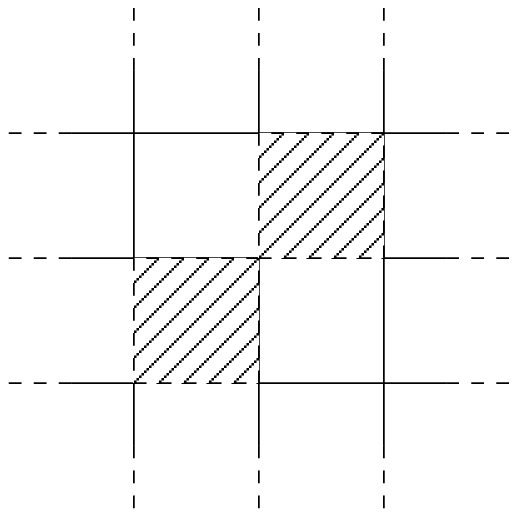} & ~~~~ &
\includegraphics[scale = 0.4]{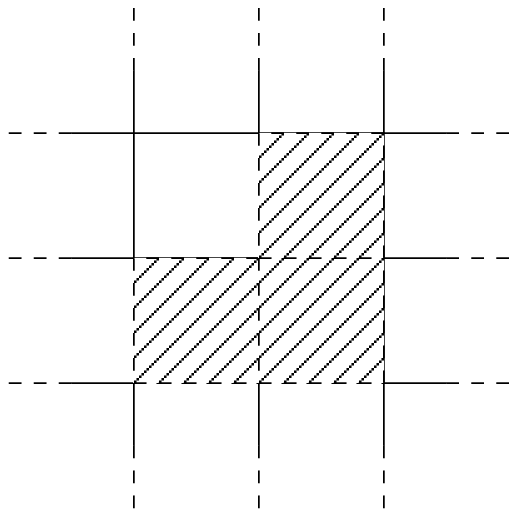}
\includegraphics[scale = 0.4]{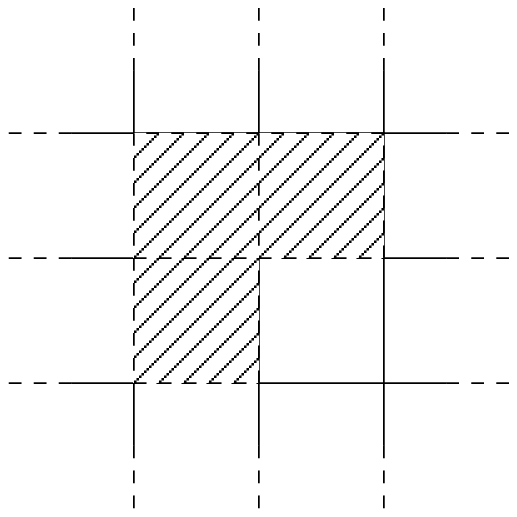} \\
(a) & ~ & (b)
\end{tabular}
\caption{(a) A diagonal step.
(b) Corners indicating that the corresponding monomials are not minimal.}
\label{fig:diacor}
\end{figure}

The shaded part in the label of a $k$-dimensional face contains the diagonal steps as in Figure \ref{fig:diacor}(a) exactly $k$ times because every time we shade in such a corner, the  
dimension decreases by one.
By straightforward counting, we get that
$$
\#~ k\text{-faces in } C_{r,n} = {r+n-k-2 \choose r-k-1, n-k-1, k}
$$
as seen in Proposition \ref{prop:fvector}.  That is, out of the $r+n-k-2$ steps we take from the  lower left corner to the upper right corner, we take $r-k-1$ steps up,  
$n-k-1$ steps right, and $k$ steps diagonally.

\begin{proposition}
The exponential generating function for the numbers $M_{r,n,k}$ of maximal $k$-faces of the tropical cyclic polytope $C_{r,n}$ is
$$
\sum_{r \geq 1, n \geq 1, k \geq 0} \frac{M_{r,n,k}}{r! n! k!} x^r y^n z^k = \frac{\partial}{\partial z} exp \left(z (y e^x - y + x e^y -x - xy)\right),
$$
and the ordinary generating function is
$$
\sum_{r \geq 1, n \geq 1, k \geq 0} M_{r,n,k} x^r y^n z^k
= \left( \frac{xy}{1-y} + \frac{yx^2}{1-x} \right) \left/\left(1- z\left(\frac{xy}{1-y} - \frac{yx^2}{1-x} \right) \right) \right. .
$$
\end{proposition}

\begin{proof}
A face is maximal if and only if the set of shaded boxes in the $r \times n$ grid does not contain any corners as in Figure \ref{fig:diacor}(b).  Then $M_{r,n,k}$ is equal to the number of $(k+1)$-tuples of either horizontal or vertical stripes of boxes, as in Figure \ref{fig:one_valent}(b), such that the sum of the widths equals $n$ and the sum of heights equals $r$.  The proposition follows from basic properties of generating functions.
\end{proof}

Moreover, every $k$-dimensional face contains precisely $2^k$ vertices  
because every diagonal step as in Figure \ref{fig:diacor}(a) 
gives $2$ ways of shading in the corners, and there are exactly $k$ such diagonal steps.  From this it is easy to see that every  
$k$-dimensional face has the combinatorial structure of a $k$-dimensional hypercube.  Therefore, the $f$-matrix of $C_{r,n}$ is very  
simple:
$f_{k,2^k} = {r+n-k-2 \choose r-k-1, n-k-1, k}$, and all other  
entries are $0$.

\section{Conclusion and Future Directions}
The methods we described can be applied to a wide range of combinatorial objects which are dual to triangulations of polytopes.  For example, tight spans of finite metric spaces can be computed using the $2 \times 2$ minors of a symmetric matrix.
There are also many enumerative questions about tropical polytopes.  For example, very little is known about the $f$-matrices.

\section{acknowledgement}
This paper grew out of a term project from the course ``Combinatorial Commuative Algebra'' taught by Bernd Sturmfels at UC Berkeley in the fall semester 2004. We are very grateful to Bernd Sturmfels for all his guidance, stimulating discussions, and inspiring questions. We also thank Mike Develin, Hiroshi Hirai, Michael Joswig, David Speyer, and the referees for helpful comments and suggestions.  Florian Block held a DAAD scholarship, and Josephine Yu was supported by an NSF Graduate Research Fellowship.

  \end{document}